\documentclass[letterpaper, 10 pt, conference]{ieeeconf}  

\IEEEoverridecommandlockouts                              
\usepackage{graphics} 
\usepackage[pdftex]{graphicx}
\usepackage{amsmath} 
\usepackage{amssymb}  
\usepackage{lilypack}

\newcommand{\R}{\mathbb{R}}
\newcommand{\T}{\text{T}}
\DeclareMathOperator*{\subjto}{subj. to}

\title{\LARGE \bf
A Dynamical Systems Approach to Energy Disaggregation
}

\author{Roy Dong, Lillian Ratliff, Henrik Ohlsson, and S. Shankar Sastry
\thanks{R. Dong, L. Ratliff, H. Ohlsson, and S. Sastry are with the Department of Electrical Engineering and Computer Sciences, University of California at Berkeley, CA, USA, {\tt\small \{roydong, ratliffl, ohlsson, sastry\}@eecs.berkeley.edu}.}
\thanks{H. Ohlsson is also with the Division of Automatic Control, Department of Electrical Engineering, Link\"oping University, Sweden.}%
\thanks{The work presented is supported by the NSF
Graduate Research Fellowship under grant DGE 1106400, NSF
CPS:Large:ActionWebs award number 0931843, TRUST (Team for Research in
Ubiquitous Secure Technology) which receives support from NSF (award
number CCF-0424422), and FORCES (Foundations Of Resilient
CybEr-physical Systems), the Swedish Research
  Council in the Linnaeus center CADICS, the European Research Council
   under the advanced grant LEARN, contract 267381, a postdoctoral grant from the Sweden-America
   Foundation, donated by ASEA's Fellowship Fund, and  by a postdoctoral
   grant from the Swedish Research Council.}%
}

\begin{document}

\maketitle
\thispagestyle{empty}
\pagestyle{empty}

\begin{abstract}

Energy disaggregation, also known as non-intrusive load monitoring (NILM), is the task of separating aggregate energy data for a whole building into the energy data for individual appliances. Studies have shown that simply providing disaggregated data to the consumer improves energy consumption behavior. However, 
placing individual sensors on every device in a home is not presently a practical solution. Disaggregation provides a feasible method for providing energy usage behavior data to the consumer which utilizes currently existing infrastructure. In this paper, we present a novel framework to perform the energy disaggregation task. We model each individual device as a single-input, single-output system, where the output is the power consumed by the device and the input is the device usage. In this framework, the task of disaggregation translates into finding inputs for each device that generates our observed power consumption. We describe an implementation of this framework, and show its results on simulated data as well as data from a small-scale experiment.

\end{abstract}

\section{INTRODUCTION}
\label{sec:introduction}
This paper is motivated by the need of efficient energy management solutions for the retail distribution domain of smart grid. Here, under the term retail distribution domain, we mean the interactions between the local distributors, e.g.\ the utility companies, and customers, e.g.\ building occupants. Usage of smart energy management devices has enabled new functionalities and has brought the potential for increased energy efficiency via real-time control and monitoring. 

Currently, we focus on commercial and residential buildings. Commercial and residential buildings are major users of energy in the developed world. Buildings account for $20$-$40$\% of total energy consumption~\cite{perez2008}. We seek to provide customers with individual device power consumption information. Studies have shown that simply providing such data improves the consumer's energy consumption behavior~\cite{Ehrhardt2010:sd}.

Current monitoring methods measure total consumption for a building. Placing individual sensors on every device in a home is not presently a practical solution. Disaggregation, also known as non-intrusive load monitoring (NILM), is the task of separating aggregate energy data for a whole building into the component energy data for individual devices, e.g.\ refrigerators, stovetops, washing machines, \&c. Disaggregation provides a feasible method for providing energy usage behavior data to the consumer, thereby allowing them to identify behavioral trends or device malfunctions that lead to inefficiencies, without requiring major infrastructural changes such as the addition of individual sensors on each device or power receptacle. 



Outside of informing consumers about ways to improve energy efficiency, disaggregation presents an opportunity for utility companies to strategically market products to consumers. It is now common practice for companies to monitor our online activity and then present advertisements which are targeted to our interests. This is known as `personalized advertising'. Disaggregation of energy data provides a means to similarly market products to consumers. This leads to the question of user privacy and the question of ownership with regards to power consumption information.
Treatment of the issue of consumer privacy in the smart grid is outside the scope of this paper. However, this is discussed in~\cite{cardenas2012}.

Additionally, disaggregation also presents opportunities for improved control. Many devices, such as heating, ventilation, and air conditioning (HVAC) units in residential and commercial buildings implement control policies that are dependent on real-time measurements. Disaggregation can provide information to controllers about system faults, such as device malfunction, which may result in inefficient control. It can also provide information about energy usage which is informative for demand response programs.

We focus on designing disaggregation methods by using 
dynamical models of devices and formulating the disaggregation problem in an optimal control framework. By working within the dynamical systems and optimal control framework, we hope that our algorithms will lend themselves to easy integration into current real-time optimal control of smart devices within the buildings and for facilitating the implementation of flexible demand response mechanisms by utilities. We designed and set up an experiment to collect data which we use for disaggregation.

The rest of the paper proceeds as follows. In Section~\ref{sec:background}, we discuss the relevant background and existing literature. In Section~\ref{sec:dynamicalmodels}, we describe our dynamical system framework for disaggregation, and implementation methods. In Section~\ref{sec:simulation}, we test our implementation on simulated data and show results. In Section~\ref{sec:experiment}, we describe the experimental setup for collecting energy data and discuss the results of the proposed disaggregation method on data from a small-scale experiment. In Section~\ref{sec:conclusion}, we make concluding remarks and discuss future work.

\section{BACKGROUND}
\label{sec:background}
The problem of non-intrusive load monitoring and the existing hardware for non-intrusive load monitoring has been studied extensively in the literature (see~\cite{berges2010:lj, berges2009:gd}). The general consensus is that non-intrusive load monitoring is a method to present the consumer with information that makes them aware of their usage and potentially provides them insight into how to improve the efficiency of their usage. Further, the technology to perform non-intrusive load monitoring is becoming widely available. Hence, there is a need for flexible and efficient disaggregation algorithms. 
  
Disaggregation of energy data has emerged as one possible solution for identifying consumer behavior patterns and device malfunctions which lead to inefficient usage of energy. The goal of the current disaggregation literature is to present methods for improving energy monitoring at the consumer level without having to place sensors at device level, but rather use existing sensors at the whole building level. 
The concept of disaggregation is not new; however, only recently has it gained attention in the energy research domain. This is likely due to the emergence of smart meters and big data analytics.

Broadly speaking, disaggregation in essence is a single-channel source separation problem. The problem of recovering the components of an aggregate signal is an inverse problem and as such is, in general, ill-posed. Most disaggregation algorithms are batch algorithms and produce an estimate of the disaggregated signals given a batch of aggregate recordings.  There have been a number of survey papers summarizing the existing methods (e.g. see \cite{Zeifman2011:fh}, \cite{kolter2011:gk}). In an effort to be as self-contained as possible, we try to provide a broad overview of the existing methods and then explain how the disaggregation method presented in this paper differs from existing solutions.  

The literature can be divided into two main approaches, namely, supervised and unsupervised. Supervised disaggregation methods require a disaggregated data set for training. This data set could be obtained by, for example, monitoring typical appliances using plug sensors. Supervised methods assume that the variations between signatures for the same type of appliances is less than that between signatures of different types of appliances. Hence, the disaggregated data set does not need to be from the building that the supervised algorithm is designed for. However, the disaggregated data set must be collected prior to deployment, and come from appliances of a similar type to those in the target building. Supervised methods are typically discriminative. 

Unsupervised methods, on the other hand, do not require a disaggregated data set to be collected. They do, however, require hand tuning of parameters, which can make it hard for the methods to be generalized in practice. It should be said that also supervised methods have tuning parameters, but these can often be tuned using the training data. 

The existing supervised methods include sparse coding \cite{kolter2010:fk}, change detection and clustering based approaches \cite{drenker1999:rk,rahayu2012:sk} and pattern recognition \cite{farinaccio1999:lk}. The sparse coding approach tries to reconstruct the aggregate signal by selecting as few signatures as possible from a library of typical signatures. Similarly, in our proposed framework we construct a library of dynamical models and reconstruct the aggregate signal by using as few as possible of these models.  

The existing unsupervised methods include factorial hidden Markov models (HMMs), difference hidden Markov models and variants \cite{kim11,Kolter:2012,Johnson2012,parson2012:AAAI:lk,pattem2012:gk} and temporal motif mining \cite{shao2012:mk}. Most unsupervised  methods models the on/off sequences of appliances using some variation of HMMs. These methods do not make use of the signature of a device and assume that the power consumption is piecewise constant. 

All method we are aware of lack the use of 
the dynamics of the devices. While the existing supervised methods often do use device signatures, these methods are discriminative and an ideal method would have a dynamical model that is capable to generating a device signature given a combination of initial state and input. Both HMMs and linear dynamical models are generative as opposed to discriminative, making them more advantageous for modeling complex system behavior. In the unsupervised domain, HMMs are used; however, they are not estimated using data and they do not model the signature of a device. The method we develop in this paper will combine the use of a generative model, i.e.\ linear dynamical models of devices, with a supervised approach to disaggregation.

\section{DYNAMICAL MODELS}
\label{sec:dynamicalmodels}
\subsection{Framework}

In our dynamical model framework, we model individual devices as single-input, single-output systems, where the output is the power consumed by the device and the input is the device usage. That is, the input is zero if a device is off, and the input is nonzero if the device is on. Thus, for device $i$, we have a model of the following form: $y_i = h_i(u_i)$, where $y_i$ is the power consumption signal of the device, $u_i$ is the input to the device, and $h_i$ is a function that represents the underlying dynamics. We build a library of models which represent the appliance types we are interested in.

With a model for each device, we treat the total power consumption as the aggregate output of all devices, i.e.\ $y = \sum_{i = 1}^D y_i$, where $y$ is the total power consumption signal and $D$ is the total number of device models. The task of disaggregation then translates into finding inputs for each device that generates our observed power consumption. In general, this solution will not be unique without more constraints on the input. Incorporating some prior on the form of the input, the problem becomes the following:
\begin{equation}
\label{eq:disagg}
\begin{array}{rl}
\arg\min_{\hat y, u} 	& L(\hat y, y_m) + g(u) 								\\
\subjto 						& \hat y_i = h_i(u_i) 								\\
								& \quad \mbox{for } i \in \{1,\dots,D\}					\\
								& \hat y = \sum_{i = 1}^D \hat y_i
\end{array}
\end{equation}
where $y_m$ is the measured power consumption, $\hat y$ is the estimated power consumption, $L$ is a loss function penalizing deviations of $\hat y$ from $y_m$, and $g$ is a regularization on the input that incorporates our priors.

\subsection{Implementation}
\label{sec:dynamic_implement}

In this framework, the task of disaggregation can be broken down into two steps: system identification and disaggregation.

\subsubsection{System identification}
In the system identification step, we seek to build a library of models which represent all the devices we are interested in. 
We assume we are given time-series measurements of power consumption for individual devices, e.g.\ a toaster, a kettle, or a LCD projector, and we wish to find a model to capture the dynamics underlying the signal. This task has a deep history and well-established literature and results~\cite{Ljung:99}.

More specifically, for some device $i$, we are given $T$ power usage samples, $y_i[k] \in \R$ for $k \in \{1,\dots,T\}$, and a sequence of corresponding inputs, $u_i[k]$ for $k \in \{1,\dots,T\}$. Assuming our world is causal, our goal is to find a satisfactory model such that $y_i[k] = h_{i,k}(u_i[1],\dots,u_i[k])$.

Throughout this paper, we will use linear time-invariant (LTI) state-space models to represent the power consumption dynamics of individual devices, i.e.\ systems of the form:
\begin{equation}
\begin{array}{rcl}
x_i[k+1] 	& = 	& A_i x_i[k] + b_i u_i[k] 	\\
y_i[k] 		& = 	& c_i^\T x_i[k]
\end{array}
\end{equation}
where $n$ is the order of the device model and $x_i[k] \in \R^n$ for $k = 1,\dots,T$ is a state underlying the dynamics. The framework generalizes to nonlinear, time-varying models as well, but for simplicity we merely consider the LTI case here.


Note that, under the assumption that similar devices have similar power consumption profiles, these models can be estimated offline. That is, for the task of disaggregation, we only need to estimate models for each class of devices once. Afterward, due to their generative nature, these models can be used for any household. 
Thus, this dynamical system framework would be cost-effective to deploy widely.

Furthermore, while power usage data can be easily recorded with plug sensors, it is not as convenient to record the input signal, $u_i[\cdot]$ for each plug. Thus, at this step, it may be necessary to apply blind system identification techniques, i.e.\ techniques for the case where both the system dynamics and the inputs are unknown. A detailed coverage of blind system identification is outside the scope of this paper; we refer the interested reader to the following references:~\cite{Abed1997,Li2006}. Also, the authors of this paper have also devised a method for blind system identification motivated by the disaggregation problem, see~\cite{Ohlssonetal:13d}.

\subsubsection{Disaggregation}
With these dynamical models in hand, we can treat disaggregation as the task of finding an input that generates our observed output. The problem formulation is as follows. We are only given samples of aggregated power consumption for a household, $y_m[k] \in \R$ for $k = 1,\dots,T$. Also, we know that the majority of the power consumption signal originates from a subset of our $D$ modeled devices. We want to find inputs which result in a similar power consumption signal.

In this paper, we take the inputs of the system to be the device's setting when it is on. Take a conventional oven as an example. It can be off, or it could be on with a temperature setting that takes on continuous values. In this situation, the input is zero if the oven is off, or the input is the temperature setting if the oven is on. An important distinction is that the input is the temperature setting, \emph{not} the temperature of the oven itself; the input can be thought of as a command to the device, e.g.\ if a user sets the oven to $350^\circ$F at time $k^*$, the input is $u_{\text{oven}}[k] = 0$ for $k < k^*$ and $u_{\text{oven}}[k] = 350$ for $k \geq k^*$. Looking at this example, we can see that a reasonable prior would be that the inputs $u_i$ are piecewise constant, and that the changes in $u_i$ across time are sparse. Throughout our paper, we use this as our prior on the inputs.

Returning to Equation \ref{eq:disagg}, we define:
\begin{equation}
	\Delta u = \begin{bmatrix}u[1]-u[0]\\u[2]-u[1]\\ \vdots \\ u[T]-u[T-1]\end{bmatrix}
\end{equation}
and we take $g(x) = \operatorname{card}(x)$, i.e.\ the number of nonzero elements in $x$.
Furthermore, we take $L$ to be the Euclidean distance on $\R^T$. Thus, we have our optimization defined.
	

A common approach when one is trying to minimize the cardinality of a vector is to relax the cardinality into the $\ell_1$ norm, which is convex. However, we found that this performs poorly in our framework. A likely explanation is that when a linear system is converted in the linear operator $\Delta u \mapsto y$, it will often fail to meet the desiderata for the $\ell_1$ relaxation.

Another technique is necessary. First, we note that if we know which elements of $\Delta u$ are nonzero, i.e.\ which devices turned on or off at what time, then it is easy to find the optimal $\Delta u$. We define each of these as a \emph{configuration}. When $g(\cdot)$ is the cardinality operator, the optimal configuration is, informally, the configuration which results in the best fit with the fewest nonzero entries. However, finding this solution is combinatorial.

We seek relaxations which will make this optimization tractable. We assume that, at each time step, only one device turns on or off at a time. This is not an egregious assumption if our sampling rate is sufficiently large. Also, we assume that the devices switch on and off in sequence; a device does not turn on and then on again afterward. We can sort by time and place our possible configurations in a tree structure. More formally, at each time step, one of $D+1$ things can happen: a device $d \in \{1,\dots,D\}$ switches on or off, where only one of the two options is possible depending on its current configuration, or no device changes configuration. This induces a hierarchical ordering on configurations of different time intervals. That is, at depth $T$ of the tree, the nodes are configurations at times $k \in \{1,\dots,T\}$, and that node's children are configurations on $\{1,\dots,T+1\}$.

If we think of the configuration at a given time as a mode, then this is a hybrid system estimation problem. The combinatorial problem above is often called a complete filter bank. This is still intractable, but we can use heuristics to intelligently prune or merge the tree and keep the set of possible configurations manageable. For the general problem, pruning and merging methods are discussed in~\cite{HwangSept.,BlomAug,Bar-ShalomMar}. These methods are known as generalized pseudo-Bayesian filters or interacting multiple models. Also, note that these algorithms allow for disaggregation to be done online.

The disaggregation problem allows for several intuitive heuristics. First, if a given configuration continues to model the future data well, we assume no device changes state. Second, if the power consumption increases by a certain amount, a device is turning on. Finally, if the power consumption decreases by a significant amount, a device is turning off. These three heuristics are sufficient to make our optimization problem extremely efficient.


\section{SIMULATION}
\label{sec:simulation}
\begin{figure}[ht]
	\begin{center}
	\includegraphics[width=\columnwidth]{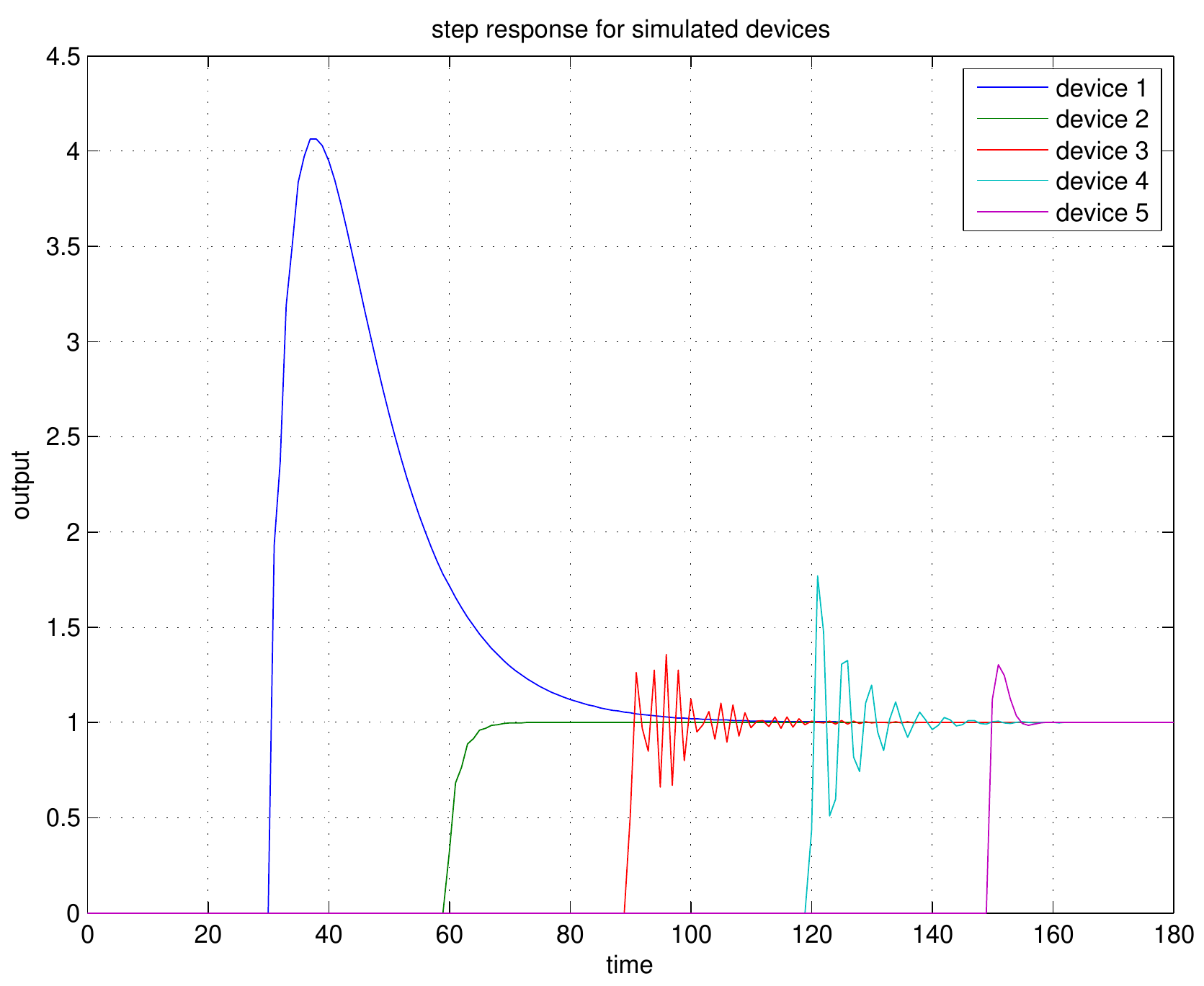}
	\end{center}
	\caption{The step responses of $D = 5$ randomly generated device models, spaced apart by 30 time steps.}
	\label{fig:sim_step}
\end{figure}

We implemented the disaggregation algorithm on simulated data. We generated $D = 5$ third-order single input, single output systems using MATLAB's {\tt drss} function, normalized to have a DC gain of 1. The step responses for these 5 systems can be seen in Figure \ref{fig:sim_step}. Let the dynamics of each system be represented with matrices $A_i, b_i, c_i^\T, d_i$ for $i \in \{1,\dots,D\}$. We assume we are given the true models for each of these $D$ devices.

\begin{figure}[ht]
	\begin{center}
	\includegraphics[width=\columnwidth]{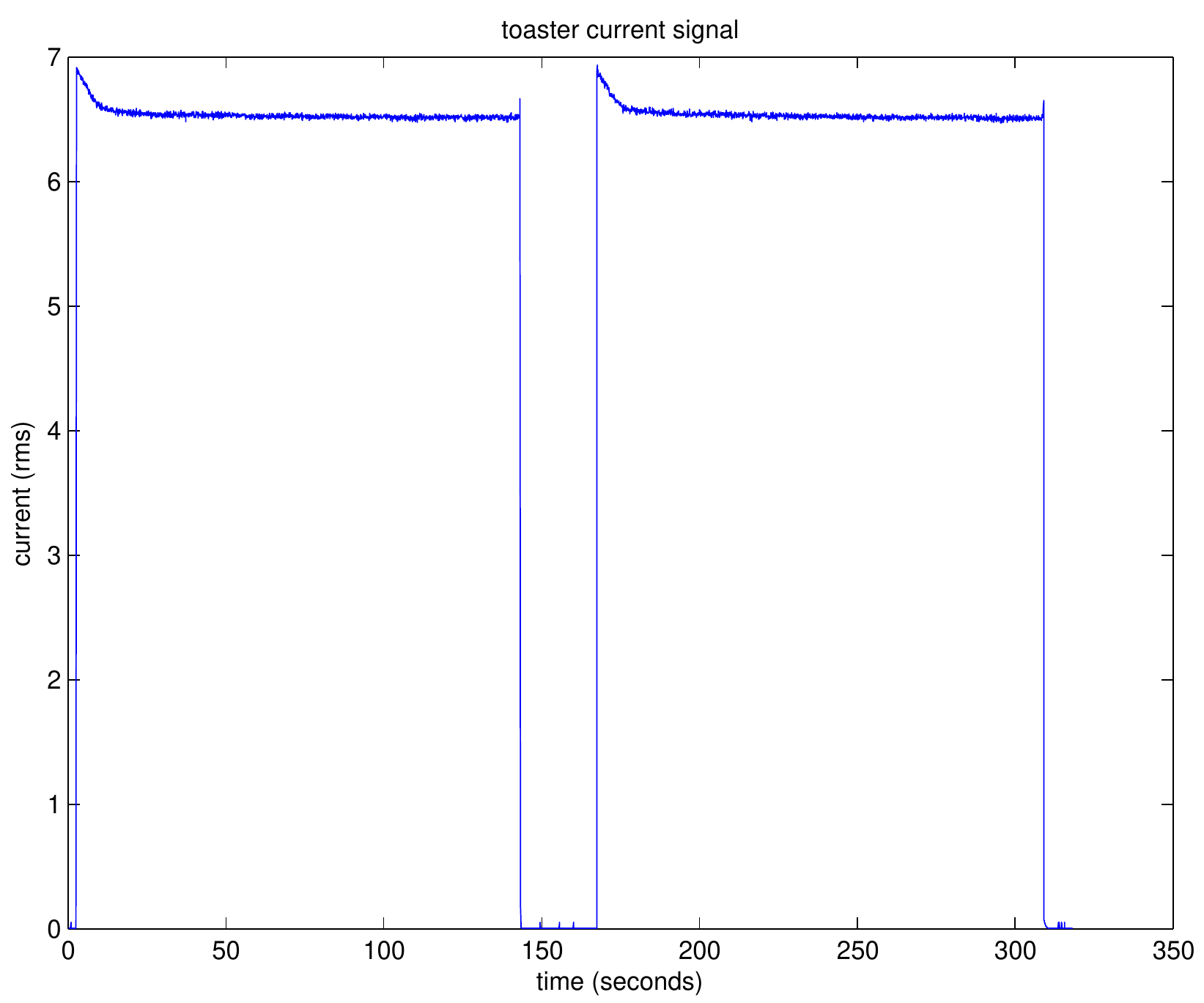}
	\end{center}
	\caption{The measured RMS current signal for a toaster. Note that the on-switches display overshoot while the off-switches do not.}
	\label{fig:toaster}
\end{figure}

We also observed that many real-life devices seem to have different dynamics when switching on and when switching off. For example, consider the root-mean-squre (RMS) current signal of a toaster, represented in Figure \ref{fig:toaster}. There is overshoot when the toaster switches on, but the off dynamics do not show the same behavior. In fact, in all of the devices we measured, we found that when a device turns off, the power drops to a negligible amount almost instantly. We factor this into our simulated models as well.

Then, we created output signals for each system by using the inputs in Equation \ref{eq:sim_inputs}. These inputs were chosen to overlap significantly. Also, not every device is activated during the simulation. The aggregated signal was created by summing these individual inputs as well as white noise with 0 mean and 0.02 standard deviation.
\begin{equation}
	\label{eq:sim_inputs}
	\begin{array}{rcl}
		u_1[k] & = & 1.2 \quad \text{for } k \in \{20,\dots,100\} \\
		u_2[k] & = & 2 \quad \text{for } k \in \{130,\dots,400\} \\
		u_3[k] & = & 0.6 \quad \text{for } k \in \{180,\dots,300\}  \\
		u_4[k] & = & 1.8  \quad \text{for } k \in \{250,\dots,350\}  \\
		u_i[k] & = & 0 \quad \text{otherwise}
	\end{array}
\end{equation}
We then run the disaggregation method described in Section \ref{sec:dynamic_implement}. For simplicity, we assume that the input is zero initially. Then, as long as this configuration's expected output and the observed output are within a certain threshold, we keep the same configuration. When the observed output exceeds this threshold, we determine if the signal is increasing or decreasing. If it is increasing, we look at all devices and nearby times to find the device that best explains the change in the measured data, as well as nearby data afterward, when driven with a constant input. If it is decreasing, since all devices turn off in the same fashion, we determine which device turned off by looking at the contribution of each device in the estimated configuration.

\begin{figure}[ht]
	\begin{center}
	\includegraphics[width=\columnwidth]{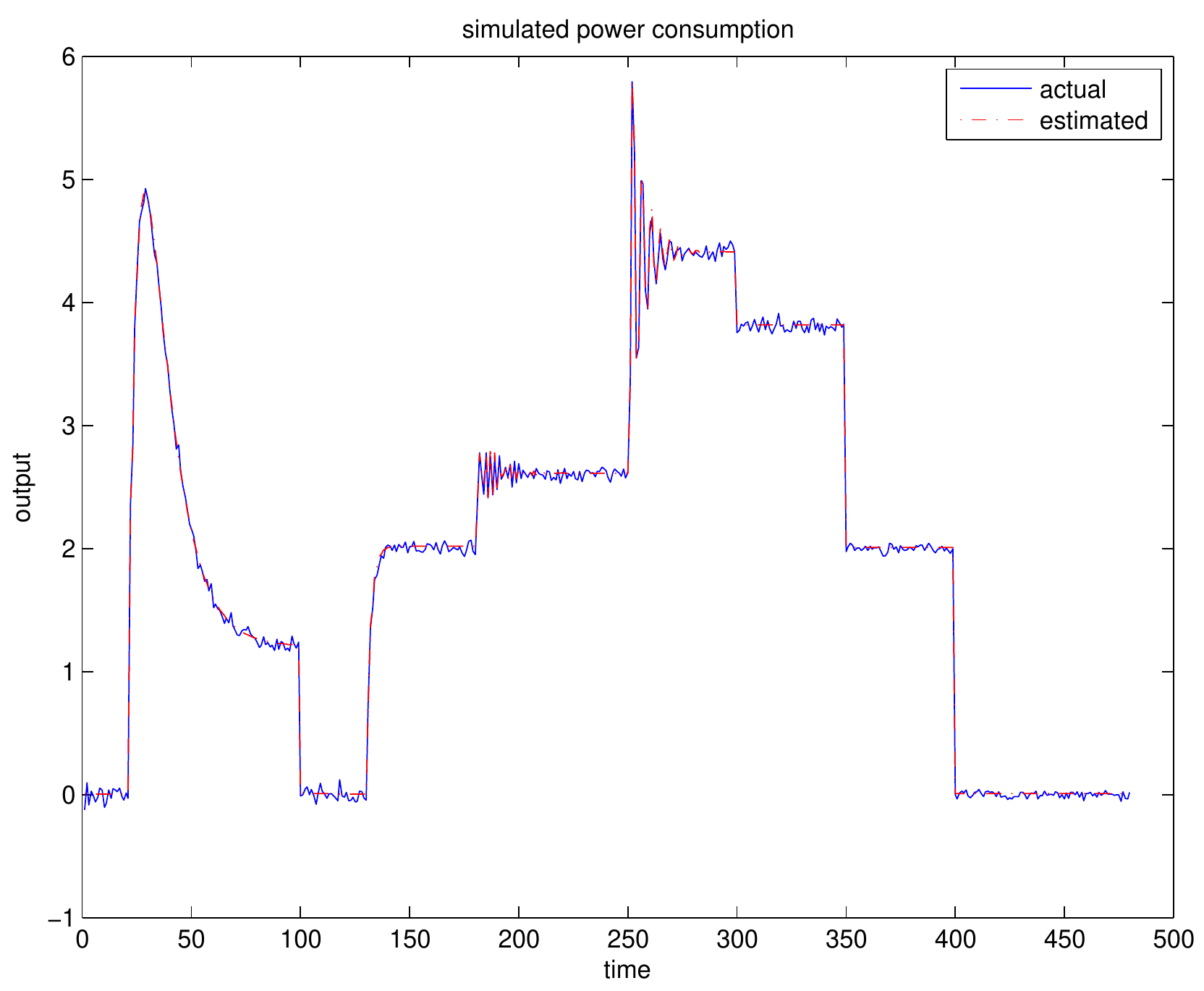}
	\end{center}
	\caption{The simulated disaggregation results.}
	\label{fig:sim_results}
\end{figure}

More formally, let $y_m$ be our measured signal, and let $\hat y$ be the predicted output under the estimated configuration. Suppose we detect a change at time $k^*$ and let $N$ be our lookahead time. Then, for nearby times $k'$ and devices $i \in \{ 1,\dots,D \}$ which are not currently on, we calculate:
\begin{equation}
\label{eq:config}
\begin{array}{rl}
\min_{y_i,x_i,u} 	& \|e - y_i\|_2^2								\\
\subjto 						& x_i[k+1] = A_i x_i[k] + b_i u 			\\
								& \quad \text{for } k \in \{k',k'+2,\dots,k^*+N-1 \} \\
								& x_i[k'] = 0 \\
								& y_i[k] = c_i^T x_i[k] + d_i u \\
								& \quad \text{for } k \in \{k',\dots,k^*+N \}
\end{array}
\end{equation}
 where $e[k] = y_m[k] - \hat y[k]$ for $k \in \{ k',\dots,k^*+N \}$ is the deviation we need to explain. 
Note here that $u$ is a scalar, not a time-dependent signal. That is, given $k'$ and $i$, we find the best input magnitude to explain the behavior. Also, note that we are implicitly reducing the cardinality of $\Delta u$, as well as reducing the number of needed calculations, by only making these estimations when our estimated configuration is not satisfactory. Furthermore, if we wish to do online disaggregation, the lookahead parameter, $N$, determines how much delay is needed.
The disaggregation estimate is:
\begin{equation}
	\begin{array}{rcl}
		\hat u_1[k] & = & 1.2017 \quad \text{for } k \in \{20,\dots,100\} \\
		\hat u_2[k] & = & 2.0104  \quad \text{for } k \in \{130,\dots,400\} \\
		\hat u_3[k] & = & 0.5827 \quad \text{for } k \in \{180,\dots,300\}  \\
		\hat u_4[k] & = & 1.7987  \quad \text{for } k \in \{250,\dots,350\}  \\
		\hat u_i[k] & = & 0  \quad \text{otherwise}
	\end{array}
\end{equation}
Every device is successfully identified, and the switching times are also correctly identified. The simulated data is plotted against the estimated data in Figure \ref{fig:sim_results}.

\section{EXPERIMENT}
\label{sec:experiment}

For the verification of our disaggregation method, we deployed a small-scale experiment.
To collect the data, we use the emonTx wireless open-source energy monitoring node from OpenEnergyMonitor\footnote{ {\tt{http://openenergymonitor.org/emon/emontx}}}. We use current transformer (CT) sensors and an alternating current (AC) to AC power adapter to measure the current and voltage respectively of the devices that we monitored. For each device we measure the root-mean-square (RMS) current ($\Irms{i}$), RMS voltage ($\Vrms{i}$), apparent power ($\Pva{i}$), real power ($\Pw{i}$), power factor ($\pf{i}$), and a UTC time stamp where the superscript $i$ index denotes the $i$th device. The sampling rate is $12$Hz.

For our experiment, we focused on small devices that would be featured in a residential or commercial office building. First, we took individual plug-level measurements for a kettle, a toaster, a projector, a monitor, and a microwave. These devices consume anywhere from $70\W$ to $1800\W$. We labeled the devices $\{ 1, \dots, 5 \}$, respectively. For the blind system identification of each of these devices, we used a simple change detection algorithm to generate input signals. Then, we fit autoregressive models with exogenous inputs. 

\begin{figure}[ht]
	\begin{center}
	\includegraphics[width=\columnwidth]{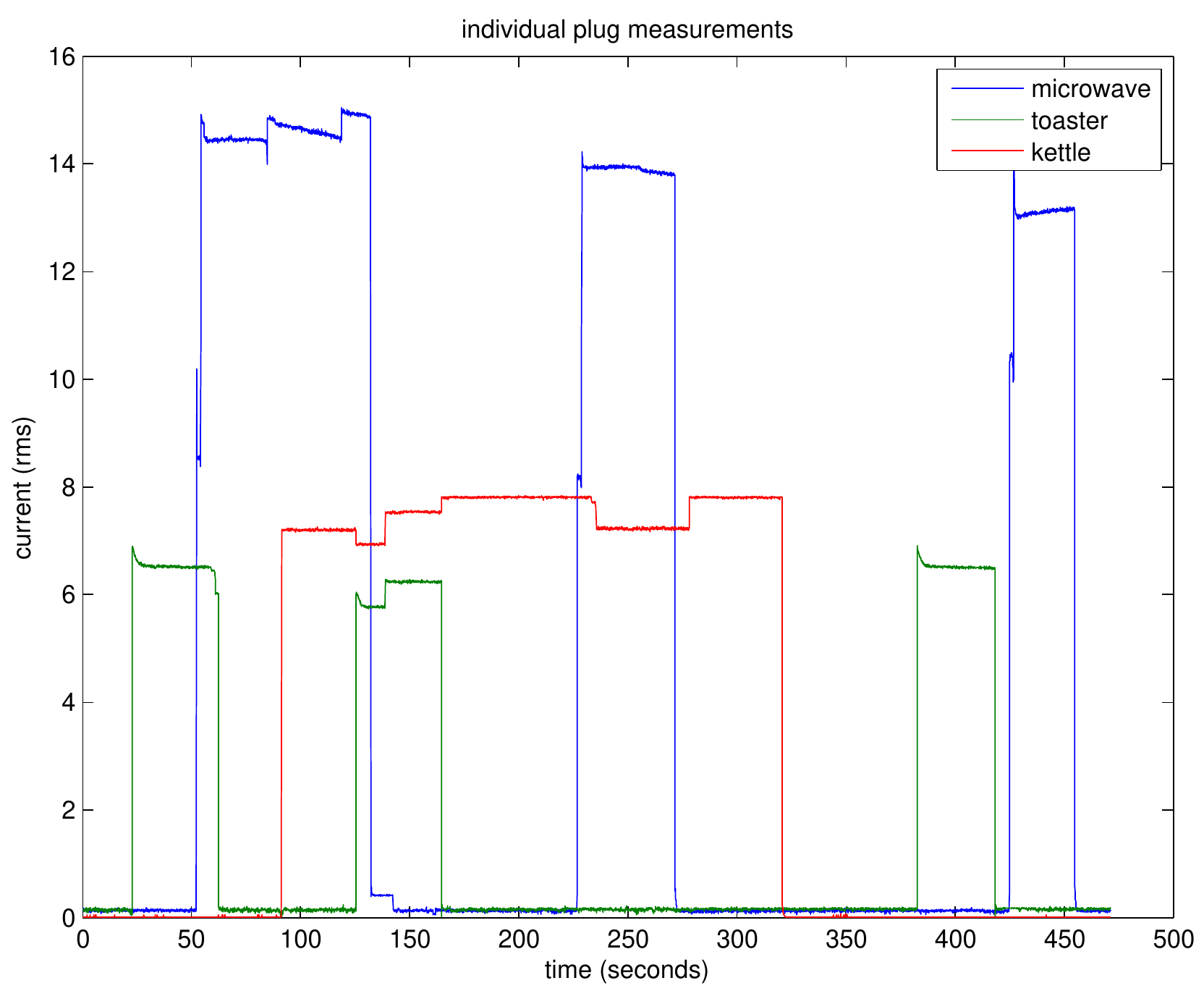}
	\end{center}
	\caption{The measurements of individual plug RMS currents.}
	\label{fig:meas_exp}
\end{figure}

Then, we ran an experiment where we had a microwave, a toaster, and a kettle (devices 5, 2, and 1, respectively) operating at different time intervals. These individual plug measurements are in Figure \ref{fig:meas_exp}. We can note that the device power consumptions are not completely independent; one device turning on can affect the power consumption of another device. However, we found this effect to be negligible in our disaggregation algorithms.


Then, we ran an experiment where we had a microwave, a toaster, and a kettle (devices 5, 2, and 1, respectively) operating at different time intervals. These individual plug measurements are in Figure \ref{fig:meas_exp}. We can note that the device power consumptions are not completely independent; one device turning on can affect the power consumption of another device. However, we found this effect to be negligible in our disaggregation algorithms.

\begin{figure}[ht]
	\begin{center}
	\includegraphics[width=\columnwidth]{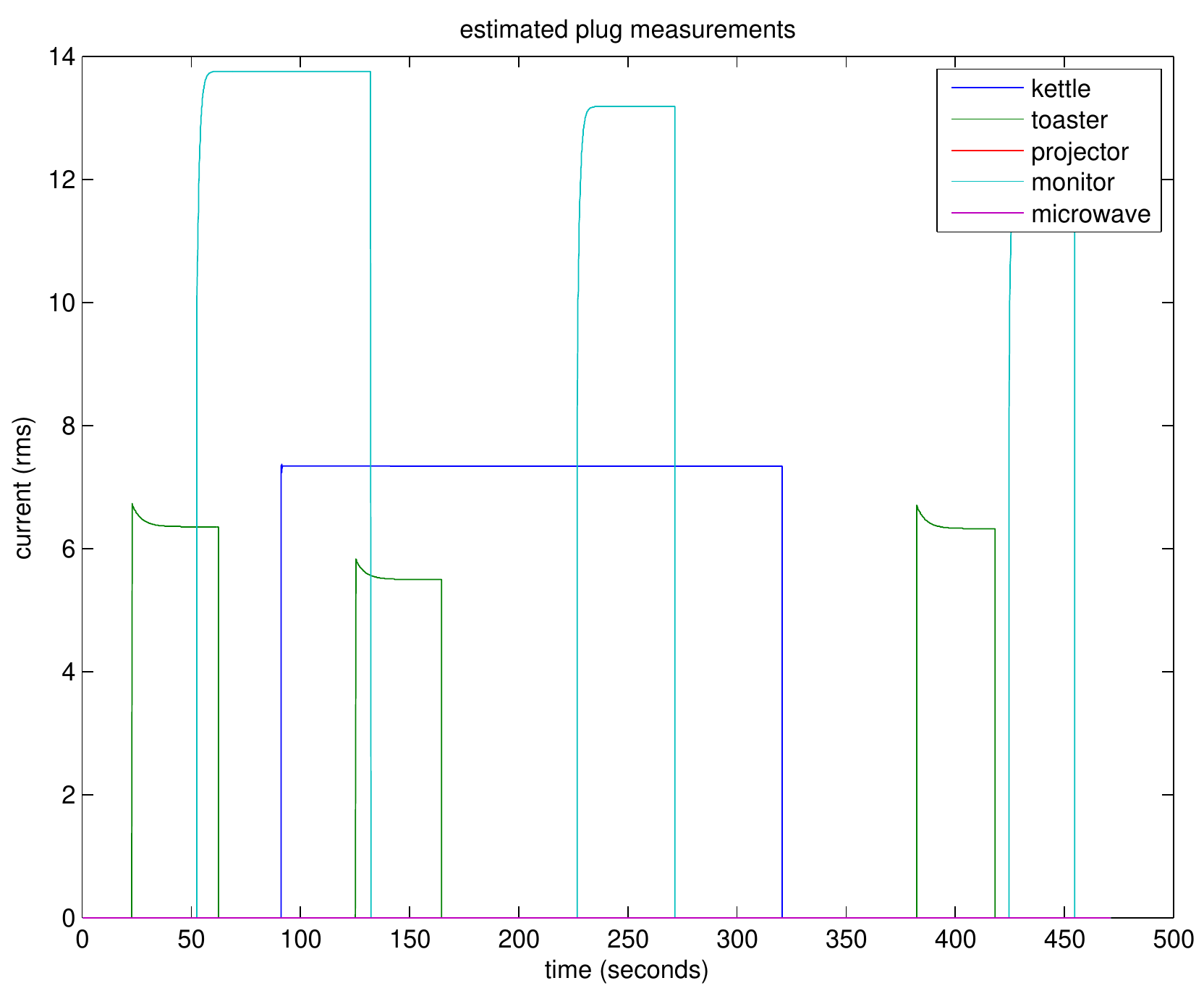}
	\end{center}
	\caption{The estimated power consumption signals of each device.}
	\label{fig:exp_inputs}
\end{figure}

The results from using our disaggregation method on the experimental data is presented in Figure \ref{fig:exp_inputs}. The estimated power consumption lines up with the measured power consumption quite well. Furthermore, the power consumption of the toaster and the kettle are correctly identified. However, the microwave is erroneously identified as a monitor. This is because the dynamics of these two models are quite similar. This error can easily be compensated for by setting a maximum power consumption for each device. That is, we can state \emph{a priori} that we know an LCD monitor will not draw over 10 amps of RMS current. When we add this prior, the microwave becomes correctly labeled.

Examining the data, we can see that methods which do not take into account the dynamics of the devices, such as the hidden Markov models (HMM) methods in~\cite{Kolter:2012,Johnson2012}, will likely confuse the kettle and toaster, which have similar amplitudes and can have similar durations. Also, the sparse coding method in~\cite{kolter2010:fk} requires a large training data set to serve as a dictionary; here, we have a very small training set from which we derive system models. Thus, a direct comparison between our method and the sparse coding method is not possible. 

\section{CONCLUSIONS AND FUTURE WORK}
\label{sec:conclusion}
In this paper, we present a novel framework to perform the task of disaggregation. We treat individual devices as systems and try to find the inputs which create the observed aggregated signal. 
This framework differs largely from the current disaggregation literature, which focuses largely on unsupervised methods. In contrast, our framework leverages many techniques and methods in system identification, optimal control, and hybrid system estimation.

We firmly believe that accounting for the power consumption profiles of individual devices will significantly improve disaggregation results. In a unsupervised setting, creating such models is very difficult. However, under the assumption that similar devices have similar power consumption profiles, the cost of collecting data and estimating these models is not significant. Thus, our framework, which utilizes more data than completely unsupervised methods, would not be infeasible to implement widely.

We tested an implementation of our framework on simulated data, as well as data from a small-scale experiment. The simulated data is completely recovered, and the experimental results closely matched the ground truth, although we did not achieve exact recovery. However, adding some reasonable assumptions allowed us to completely recover the ground truth.

For future work, we plan on deploying an experiment where we collect measurements from more devices in an actual residential setting. In this experimental setting, we hope to learn not only device dynamics, but also the user's consumption patterns. One of the benefits of our framework is that we can learn devices independent of the consumer, and then learn the user's consumption patterns. Note that in many unsupervised methods, keeping the device constant while varying the consumer's usage patterns would result in different models entirely.

Additionally, throughout our experiments, we noticed that some devices do not fit our current modeling assumptions neatly. For example, the microwave warms up for a second or two, and begins heating. This results in two successive jumps in power consumption. With our current modeling assumptions, the best fit is an over-damped system. This is not ideal, and we hope to model devices as hybrid systems with multiple modes in the future.

\section{ACKNOWLEDGMENTS}
The authors would like to thank Aaron Bestick for his advice and many helpful discussions.


\bibliographystyle{IEEEtran}
\bibliography{main}

\end{document}